\pdfoutput=1

\documentclass[twoside,12pt,a4paper]{amsart} 
\usepackage{amsmath,amsfonts,amsthm,eurosym}
\usepackage[colorlinks=true,linkcolor=blue,citecolor=blue,filecolor=blue,urlcolor=blue,pageanchor=true,plainpages=false,linktocpage]{hyperref}
\usepackage[english]{babel}
\usepackage{setspace}
\usepackage{xcolor}

\numberwithin{equation}{section}

\renewcommand{\epsilon}{\varepsilon}
\renewcommand{\phi}{\varphi}
\renewcommand{\rho}{\varrho}
\renewcommand{\theta}{\vartheta}

\begin{document}
\baselineskip3.6ex

\title{}

\author{}

\makeatletter
\def\@makefnmark{}
\makeatother
\newcommand{\myfootnote}[2]{\footnote{\textbf{#1}: #2}}
 \footnote {\textit{2020 Mathematics Subject Classification:} 53C25, 53C21 \\ \textit{Keywords:} Semi-Riemannian metrics, Einstein manifolds, sequential warped-product of special type, positive constant Ricci curvature. \\  \\ \textit{$[1]$}: Mathematical and Physical Science Foundation, Sidevej 5, 4200 Slagelse, Denmark.\\ $[2]$: King Abdulaziz University, Department of Mathematics, 21589 Jeddah KSA. \\ $[3]$: Mathematical and Physical Science Foundation, Sidevej 5, 4200 Slagelse, Denmark. \\ $[4]$: Institute of Experimental Physics, Slovak Academy of Sciences, Kosice, Slovak Republic. \\ $[5]$: Simon Fraser University, Burnaby, B.C., Canada, V5A 1S6 , and, The Pacific Institute for the Mathematical Sciences.}
\title{AA}
{\bfseries\centerline{A family of special case of sequential warped-product manifolds}}
{\bfseries\centerline{with semi-Riemannian Einstein metrics}}
\\
\\
\centerline{Alexander Pigazzini$^{[1]}$, Cenap \"{O}zel$^{[2]}$, Saeid Jafari$^{[3]}$,}
\\
\centerline{Richard Pincak$^{[4]}$ and Andrew DeBenedictis$^{[5]}$}
\\
\\
\\{\bfseries \centerline {Abstract}} 
\\ 
\\
\\
We derive the general formulas for a special configuration of the sequential warped-product semi-Riemannian manifold to be Einstein, where the base-manifold is the product of two manifolds both equipped with a generic diagonal conformal metrics. Subsequently we study the case in which these two manifolds are conformal to a $n_1$-dimensional and $n_2$-dimensional pseudo-Euclidean space, respectively. For the latter case,  we prove the existence of a family of solutions that are invariant under the action of a $(n_1-1)$-dimensional group of transformations to the case of positive constant Ricci curvature ($\lambda>0$).
\\
\\
\\
\\
{\bfseries \centerline {1. Introduction and Preliminaries}}
\\
\\
The warped-product manifolds are type of manifolds introduced by Bishop and O'Neill \cite{Bishop}.
These manifolds have become very important in the context of differential geometry and are also extensively studied in the arena of General Relativity, for instance with respect to generalized Friedmann-Robrtson-Walker spacetimes. Many properties for warped product manifolds and submanifolds were presented by B.-Y. Chen in \cite{Chen}.
\\
A warped-product manifold can be constructed as follows. Let $(B, g_B)$ and $(F, g_F)$ be two semi-Riemannian manifolds and $\tau$, $\sigma$ be the projection of $B \times F$ onto $B$ and $F$, respectively.
\\
The warped-product $M=B \times_f F$ is the manifold $B \times F$ equipped with the metric tensor $g=\tau^*g_B+f^2 \sigma^*g_F$, where $^*$ denotes the pullback and $f$ is a positive smooth function on $B$, the so-called warping function.
\\
\\
Explicitly, if $X$ is tangent to $B \times F$ at $(p,q)$ (where $p$ is a point on $B$ and $q$ is a point on $F$), then:
\\
{\centerline {$\langle X,X \rangle=\langle d\tau(X),d\tau(X)\rangle+f^2(p)(d\sigma(X),d\sigma(X))$.}}
\\
$B$ is called the \textit{base-manifold} of $M=B \times_f F$ and $F$ is the \textit{fiber-manifold}. If $f=1$, then $B \times_f F$ reduces to a semi-Riemannian product manifold. The leaves $B \times q = \sigma^{-1}(q)$ and the fibers $p \times F =\tau^{-1}(p)$ are Riemannian submanifolds of $M$. Vectors tangent to leaves are called horizontal and those tangent to fibers are called vertical. $By \; \mathcal{H}$ we denote the orthogonal projection of $T_{(p,q)}M$ onto its horizontal subspace $T_{(p,q)}(B \times q)$ and $\mathcal{V}$ denotes the projection onto the vertical subspace $T_{(p,q)}(p \times F)$, see \cite{O'Neill}.
\\
\\
If $M$ is an $n$-dimensional manifold, and $g_M$ is its metric tensor, the Einstein condition means that $Ric_M = \lambda g_M$ for some constant $\lambda$, where $Ric_M$ denotes the Ricci tensor of $g_M$. An Einstein manifold with $\lambda = 0$ is called Ricci-flat manifolds.
\\
Then keeping this in mind, we get that a warped-product manifold $(M, g_M)=(B,g_B)\times_f(F,g_F)$ (where ($B, g_B$) is the base-manifold, ($F, g_F$) is the fiber-manifold), with 
\\
$g_M=g_B+f^2 g_F$,  is Einstein if only if (see \cite{Chen}):
\\
\\
\numberwithin{equation}{section}
{(1.1)}
$Ric_M=\lambda g_M \Longleftrightarrow\begin{cases} 
 Ric_B- \frac{d}{f}Hess(f)= \lambda g_B  \\  Ric_F=\mu g_F \\ f \Delta f+(d-1) |\nabla f|^2 + \lambda f^2 =\mu
\end{cases}$
\\
\\
where $\lambda$ and $\mu$ are constants, $d$ is the dimension of $F$, $Hess(f)$, $ \Delta f$ and $\nabla f$ are, 
\\
respectively, the Hessian, the Laplacian (given by $tr \; Hess(f)$) and the gradient of $f$ for $g_B$, with $f:(B) \rightarrow \mathbb{R}^+$ a smooth positive function.
\\
\\
Contracting first equation of (1) we get: 
\\
\\
\numberwithin{equation}{section}
{(1.2)}
$R_Bf^2-f \Delta fd=n f^2 \lambda$ 
\\
where $n$ and $R_B$ is the dimension and the scalar curvature of $B$ respectively. From third equation, considering $d \neq 0$ and $d \neq 1$, we have:
\\
\\
\numberwithin{equation}{section}
{(1.3)}
$f\Delta fd+d(d-1)|\nabla f|^2+\lambda f^2d=\mu d$
\\
Now from (1.2) and (1.3) we obtain:
\\
\numberwithin{equation}{section}
{(1.4)}
$|\nabla f|^2+[\frac{\lambda (d-n)+R_B}{d(d-1)}]f^2=\frac{\mu}{(d-1)}$.
\\
\\
In 2017 de Sousa and Pina \cite{Sousa}, studied warped-product semi-Riemannian Einstein manifolds in case that base-manifold is conformal to an n-dimensional pseudo-Euclidean space and invariant under the action of an $(n-1)$-dimensional group with Ricci-flat fiber $F$. In \cite{Pal} the authors extend the work done for multiply warped space.
\\
\\
In \cite{Shen}, the author introduced a new type of warped-products called sequential warped-products, i.e. $(M, g_M)$ where $M=(B_1 \times_h B_2) \times_{f}F$ and $g_M=(g_{B_1} + h^2g_{B_2})+f^2g_F$, to cover a wider variety of exact solutions to Einstein’s field equation.
\\
Regarding the sequential warped-product manifolds, some works have been published in recent years (\cite{Chand Dea}, \cite{Sahin}, \cite{Guler}, \cite{Pahan}, \cite{Karaca}, \cite{Kumar}).
\\
\\
The main aim of the present paper is largely to continue to extend the work done in \cite{Sousa} (as was done for the multiply warped-product manifold in \cite{Pal}), also for a special case of sequential warped-product manifolds, (i.e. for $h=1$, with $B_2$ as an Einstein manifold, and flat fiber $F$, where the base-manifold $B=B_1\times B_2$ is the product of two manifolds both equipped with a conformal metrics, and the warping function is a smooth positive function $f(x,y)=f_1(x)+ f_2(y)$ where each is a function on its individual manifold). The method will be as follows: first deriving the general formulas to be Einstein and second, providing the existence of solutions that are invariant under the action of a $(n_1-1)$-dimensional group of transformations to the case of positive constant Ricci curvature. In fact, since in both references, \cite{Sousa} and \cite{Pal}, the authors show solutions for the Ricci-flat case ($\lambda = 0$), we, following their same construction, show the existence of a family solutions for constant positive Ricci curvature ($\lambda> 0$).  In particular, this proof of the existence of a family of solutions also holds for \cite{Sousa} considering $dimF = dimB$.
\\
\\
{\bfseries Definition 1.1:} We consider the special case of the Einstein sequential warped-product manifold, that satisfies (1.1). The manifold $(M, g_M)$ comprises the base-manifold $(B, g_B)$ which is a Riemannian (or pseudo-Riemannian) product-manifold $B=B_1\times B_2$, with $B_2$ as an Einstein manifold (i.e., $Ric_{B_2}=\lambda g_{B_2}$, where $\lambda$ is the same for (1.1) and $g_{B_2}$ is the metric for $B_2$), and $dim(B_2)=n_2$, $dim(B_1)=n_1$ the dimension of $B_2$ and $B_1$, respectively, so that $dim(B)=n=n_1+n_2$. The warping function $f:B \rightarrow \mathbb{R}^+$ is a smooth positive function $f(x,y)=f_1(x)+ f_2(y)$  (where each is a function on its individual manifold, i.e., $f_1:B_1\rightarrow \mathbb{R}^+$ and $f_2:B_2 \rightarrow \mathbb{R}^+$). The fiber-manifold $(F, g_F)$ is the $\mathbb{R}^d$, with orthogonal Cartesian coordinates such that $g_{ab}=-\delta_{ab}$.
\\
\\
{\bfseries Proposition 1.1:} If we write the B-product as $B=B_1\times B_2$, where:
\\
i) $Ric_{B_i}$ is the Ricci tensor of $B_i$ referred to $g_{B_i}$, where $i=1, 2$,  
\\
ii) $f(x, y)=f_1(x)+ f_2(y)$, is the smooth warping function, where $f_i:B_i\rightarrow \mathbb{R}^+$,
\\
iii) $Hess(f)=\sum_i {\tau_i^*Hess_i(f_i)}$ is the Hessian referred on its individual metric, where $\tau_i^*$ are the respective pullbacks, (and $\tau_2^* Hess_2(f_2)=0$ since $B_2$ is Einstein), 
\\
iv) $\nabla f$ is the gradient (then $|\nabla f|^2= \sum_i{|\nabla_if_i|^2}$), and
\\
v) $\Delta f=\sum_i{\Delta_if_i}$ is the Laplacian, (from (iii) therefore also $\Delta_2 f_2=0$). 
\\
Then the Ricci curvature tensor will be:
\\
\\
\numberwithin{equation}{section}
{(1.5)}
$\begin{cases} 
Ric_{M}(X_i, X_j) = Ric_{B_1}(X_i, X_j) - \frac{d}{f}Hess_1(f_1)(X_i, X_j)\\  Ric_{M}(Y_i, Y_j) = Ric_{B_2}(Y_i, Y_j)  \\ Ric_{M}( U_i, U_j)=Ric_{F}( U_i, U_j) -g_F( U_i, U_j)f^*\\  Ric_{M}(X_i, Y_j)=0 \\ Ric_{M}(X_i, U_j)=0, \\ Ric_{M}(Y_i, U_j)=0,
\end{cases}$
\\
where  $f^*= \frac{\Delta_1 f_1}{f}+(d-1) \frac{|\nabla f|^2}{ f^2}$, and $X_i$, $X_j$, $Y_i$, $Y_j$, $U_i$, $U_j$ are vector fields on $B_1$, $B_2$ and $F$, respectively.
\\
\\
{\bfseries Theorem 1.1:} \textit{ A warped-product manifold is a special case of an Einstein sequential warped-product manifold, as defined in \textit{Definition 1.1}, if and only if:
\\
\\
\numberwithin{equation}{section}
{(1.6)}
$Ric_{M}=\lambda g_M \Longleftrightarrow\begin{cases} 
 Ric_{B_1}- \frac{d}{f}\tau_1^*Hess_1(f_1)= \lambda g_{B_1}\\ \tau_2^* Hess_2(f_2)=0
\\
Ric_{B_2} = \lambda g_{B_2} \\ Ric_{F}=0 \\ f \Delta_1 f_1+(d-1) |\nabla f|^2 + \lambda f^2 =0,
\end{cases}$
\\
(since $Ric_B$ is the Ricci curvature of $B$ referred to $g_B$, then $Ric_B=Ric_{B_1}+Ric_{B_2}=\lambda(g_{B_1}+g_{ B_2})+ \frac{d}{f}\tau_1^*Hess_1(f_1)$}.
\\
\\
Therefore from (1.2) and (1.3):
\\
\\
\numberwithin{equation}{section}
{(1.7)}
$R_{M}=\lambda (n+d) \Longleftrightarrow\begin{cases} 
R_{B_1}f- \Delta_1f_1d=n_1 f \lambda \\ \Delta_2 f_2=0 \\ R_{B_2}= \lambda n_2 \\ R_{F}=0 \\ f \Delta_1f_1+(d-1) |\nabla f|^2 + \lambda f^2 =0.
\end{cases}$
\\
\\
where $n_1$ and $R_1$ are the dimension and the scalar curvature of $B_1$ referred to $g_{B_1}$, respectively.
\\
\\
\\
\textit{Proof.} We applied the condition that the warped-product manifold of system (1.5) is Einstein.$\qed$
\\
\\
This particular type of Einstein sequential warped-product manifold, as per \textit{Definition 1.1}, allows to cover a wider variety of exact solutions of Einstein's field equation, without complicating the calculations much, compared to the Einstein warped-product manifolds with Ricci-flat fiber ($F, g_F$), also considered by the authors of \cite{Sousa}.
\\
\\
{\bfseries \centerline {2. Conformal B-metrics}}
\\
\\
In this section we will consider a special type of sequential warped-product manifold $(M, g_M)$, as described in the previous section, but in which the base-manifold is the product of two manifolds both equipped with a conformal metrics. First we will show the general formulas for which such a manifold $M$ is Einstein, then we will show the same in the case where the conformal metrics are both diagonal, and finally for the case in which the base-manifold is the product of two conformal manifolds to a $n_1$-dimensional and $n_2$-dimensional pseudo-Euclidean space, respectively.
\\
\\
{\bfseries Theorem 2.1:} \textit{ Let $(B, g_B)$, be the base-manifold $B=( B_1 \times B_2)$, $B_1=\mathbb{R}^{n_1}$, with coordinates $(x_1, x_2, .. x_{n_1})$, $B_2=\mathbb{R}^{n_2}$, with coordinates $(y_1, y_2, .. y_{n_2})$, where  $n_1, n_2 \geq 3$, and let  $g_B= g_{B_1} + g_{B_2}$ be the metrics on $B$, where $g_{B_1}=\epsilon_i \delta_{ij}$ and $g_{B_2}=\epsilon_l \delta_{lr}$. 
\\
Let $f_1:\mathbb{R}^{n_1} \rightarrow \mathbb{R}$, $f_2: \mathbb{R}^{n_2} \rightarrow \mathbb{R}$, $\phi_1:\mathbb{R}^{n_1} \rightarrow \mathbb{R}$ and ${\phi_2}:\mathbb{R}^{n_2} \rightarrow \mathbb{R}$, be smooth functions, where $f_1$ and $f_2$ are positive functions, such that $f=f_1+ f_2$ as in \textit{Definition 1.1}. Finally, let $(M, g_M)$ be $((B_1 \times B_2)\times_{f=f_1 + f_2}F, g_M)$, with $g_M=\bar g_B+(f_1+f_2)^2 g_F$, with conformal metric $\bar g_B=\bar g_{B_1} + \bar g_{B_2}$, where $\bar g_{B_1}=\frac{1}{\phi_1^2}g_{B_1}$, $\bar g_{B_2}=\frac{1}{{\phi_2}^2}g_{B_2}$, and $F=\mathbb{R}^{d}$ with $g_F=-\delta_{ab}$.
\\
Then the warped-product metric $g_M=\bar g_B+(f_1+f_2)^2 g_F$ is Einstein with constant Ricci curvature $\lambda$ if and only if, the functions $f_1$, $f_2$, $\phi_1$ and $\phi_2$ satisfy:
\\
\\
\textit{(I)} $(n_1-2)f\phi_{1_{, x_i x_j}}-\phi_1f_{1_{, x_i x_j}}d - \phi_{1_{, x_i}}f_{1_{, x_j}}d-\phi_{1_{, x_j}}f_{1_{, x_i}}d=0$ for $i \neq j$,
\\
\\
\textit{(II)} $(n_2 - 2) \phi_{2_{, y_l y_r}}=0$ for $l \neq r$,
\\
\\
\textit{(III)} $\phi_1[(n_1-2)f \phi_{1_{, x_ix_i}} - \phi_1f_{1_{, x_ix_i}}d-2\phi_{1_{, x_i}}f_{1_{, x_i}}d]+$
\\
\\
{\centerline {$+\epsilon_i [f \phi_1 \sum_{k=1}^{n_1} \epsilon_k\phi_{1_{, x_k x_k}}-(n_1-1)f \sum_{k=1}^{n_1} \epsilon_k {\phi_1}^2_{, x_k}+\phi_1 d \sum_{k=1}^{n_1}\epsilon_k \phi_{1_{, x_k}}f_{1_{, x_k}}]=\epsilon_i \lambda f$,}}
\\
\\
\textit{(IV)} $\phi_2(n_2-2)\phi_{2_{, y_l y_l}}+\epsilon_l \phi_2\sum_{s=1}^{n^2}\epsilon_s\phi_2{_{, y_s y_s}}-(n_2 -1) \epsilon_l \sum_{s=1}^{n_2}\epsilon_s {\phi_2}^2_{, y_s}=\lambda \epsilon_l$,
\\
\\
\textit{(V)} $-f {\phi_1}^2 \sum_{k=1}^{n_1} \epsilon_k f_{1_{, x_k x_k}}+(n_1-2)f \phi_1 \sum_{k=1}^{n_1}\epsilon_k \phi_{1_{, x_k}}f_{1_{, x_k}} +$
\\
\\
{\centerline {$-(d-1)({\phi_1}^2 \sum_{k=1}^{n_1} \epsilon_k {f_1}^2_{, x_k}+{\phi_2}^2 \sum_{s=1}^{n_2}\epsilon_s {f_2}^2_{, y_s})=\lambda f^2$.}}}
\\
\\
\\
Before proving \textit{Theorem 2.1}, and showing the existence of a solution for $\lambda> 0$, we want to deduce the formulas for generic diagonal conformal metrics $g_{B_1}$ and $g_{B_2}$. 
\\
Based on this, we consider $(B, g_B)$, the base-manifold $B=( B_1 \times B_2)$, with $dim(B_1)=n_1$, $dim(B_2)= n_2$, and $g_B= g_{B_1} + g_{B_2}$. We also consider $f_1:\mathbb{R}^{n_1} \rightarrow \mathbb{R}$, $f_2: \mathbb{R}^{n_2} \rightarrow \mathbb{R}$, $\phi_1:\mathbb{R}^{n_1} \rightarrow \mathbb{R}$ and ${\phi_2}:\mathbb{R}^{n_2} \rightarrow \mathbb{R}$, are smooth functions, where $f_1$ and $f_2$ are positive functions, such that $f=f_1+ f_2$ as in \textit{Definition 1.1}. And finally, we consider $(M, g_M)$ with $((B_1 \times B_2)\times_{(f_1 + f_2)}F, g_M)$, with $g_M=\bar g_B+(f_1+f_2)^2 g_F$, with conformal metric $\bar g_B=\bar g_{B_1} + \bar g_{B_2}$, where $\bar g_{B_1}=\frac{1}{\phi_1^2}g_{B_1}$, $\bar g_{B_2}=\frac{1}{{\phi_2}^2}g_{B_2}$, and $F=\mathbb{R}^{d}$ with $g_F=-\delta_{ab}$.
\\
\\
From (1.6), considering the conformal metric on $B_1$ and $B_2$, it is easy to deduce that $M$ is Einstein if and only if: 
\\
\numberwithin{equation}{section}
{(2.1)}
$Ric_{\bar B_1}=\lambda \bar g_{B_1} + \frac{d}{f}Hess_{\bar1}(f _1)$, or equivalently \numberwithin{equation}{section}
{(2.2)}
$R_{\bar B_1}=\lambda n_1 + \frac{d}{f}\Delta_{\bar1}(f_1)$,
\\
\numberwithin{equation}{section}
{(2.3)}
$Ric_{\bar B_2}=\lambda \bar g_{B_2}$, or equivalently \numberwithin{equation}{section}
{(2.4)}
$R_{\bar B_2}=\lambda n_2$,
\\
\numberwithin{equation}{section}
{(2.5)}
$0= \lambda f^2 + f \Delta_{\bar 1}f_1+(d-1)[|\nabla_{\bar 1}f_1|^2 + | \nabla_{\bar 2}f_2|^2]$.
\\
\\
If we consider a generic diagonal metric, $\bar g_{B_{ij}}=\bar g_{B_{1_{ij}}}+\bar g_{B_{2{ij}}}=\eta_{ij}$, and $\eta_{ij} =0$ for $i \neq j$, then $M$ is Einstein if and only if (2.1), (2.3) (or equivalently (2.2), (2.4)), (2.5) and the following, are satisfied:
\\
{(2.6)}
$Ric_{\bar B_1}= \frac{d}{f}Hess_{\bar 1}(f _1)$, for $i \neq j$,
\\
{(2.7)}
$Ric_{\bar B_2}=0$, for $i \neq j$.
\\
\\
\textit{Proof of Theorem 2.1.} At this point we can calculate:
\\
\numberwithin{equation}{section}
{(2.8)}
$Ric_{\bar B_1}=\frac{1}{\phi_1^2}\{(n_1-2)\phi_1 Hess_1(\phi_1) +[\phi_1 \Delta_1 \phi_1-(n_1-1)|\nabla_1\phi_1|^2]g_{B_1}\}$,
\\
\numberwithin{equation}{section}
{(2.9)}
$Ric_{\bar B_2}=\frac{1}{\phi_2^2}\{(n_2-2) \phi_2 Hess_2(\phi_2) +[\phi_2 \Delta_2 \phi_2-(n_2-1)| \nabla_2 \phi_2 ^2]g_{B_2}\}$,
\\
so we can write:
\\
\numberwithin{equation}{section}
{(2.10)}
$Ric_{\bar B_1}(X_i, X_j)=\frac{1}{\phi_1^2}\{(n_1-2)\phi_1 Hess_1(\phi_1)(X_i, X_j) +[\phi_1 \Delta_1\phi_1-(n_1-1)|\nabla_1\phi_1|^2]g_{B_1}(X_i, X_j)\}$,
\\
\numberwithin{equation}{section}
{(2.11)}
$Ric_{\bar B_2}(Y_l, Y_r)=\frac{1}{\phi_2^2}\{(n_2-2) \phi_2 Hess_2 (\phi_2)(Y_l, Y_r) +[\phi_2 \Delta_2 \phi_2-(n_2-1)| \nabla_2 \phi_2|^2]g_{B_2}(Y_l, Y_r)\}$,
\\
\numberwithin{equation}{section}
{(2.12)}
$Ric_{M}(X_i, X_j)=Ric_{\bar B_1}(X_i, X_j) - \frac{d}{f} Hess_{\bar 1}(f_1)(X_i, X_j)$,
\\
for what was stated in \textit{Proposition 1.1} we have:
\\
\numberwithin{equation}{section}
{(2.13)}
$Ric_{M}(Y_l, Y_r)=Ric_{\bar B_2}(Y_l, Y_r)$,
\\
and in the end
\\
\numberwithin{equation}{section}
{(2.14)}
$Ric_{M}(X_i, Y_j)=0$.
\\
\numberwithin{equation}{section}
{(2.15)}
$Ric_{M}(X_i, U_j)=0$.
\\
\numberwithin{equation}{section}
{(2.16)}
$Ric_{M}(Y_i, U_j)=0$.
\\
Since $Ric_{F}=0$ we obtain:
\\
\numberwithin{equation}{section}
{(2.17)}
$Ric_{M}(U_i, U_j)= -g_M(U_i, U_j)(\frac{\Delta _{\bar 1}f_1}{f}+(d-1)\frac{g_M(\nabla f, \nabla f)}{f^2})$,
\\
where, analogous to \textit{Proposition 1.1}, we consider $g_M(\nabla f, \nabla f)= \bar g_{B_1}(\nabla f_1, \nabla  f_1)+ \bar g_{B_2}(\nabla f_2, \nabla f_2)$.
\\
\\
 Let $\phi_{1_{, x_i x_j}}$, $\phi_{1_{, x_i}}$, $f_{1_{, x_i x_j}}$, $f_{1_{, x_i}}$, $\phi_{2_{, y_l y_r}}$, $\phi_{2_{, y_l}}$, $f_{2_{, y_l y_r}}$ and $f_{2_{, y_l}}$, be the second and the first order derivatives of $\phi_1$, $\phi_2$, $f_1$ and $f_2$,  respectively, with respect to $x_ix_j$ and $y_ly_r$.
\\
Now we have:
\\
{(2.18)}
$Hess_{1}(\phi_1)(X_i, X_j)=\phi_{1_{, x_ix_j}}$, 
\\
{(2.19)}
$\Delta_{1}(\phi_1)=\sum_{k=1}^{n_1} \epsilon_k \phi_{1_{, x_k x_k}}$, 
\\
{(2.20)}
$|\nabla_{1}(\phi_1)|^2=\sum_{k=1}^{n_1} \epsilon_k \phi_{1_{, x_k}}^2$, 
\\
{(2.21)}
$Hess_{2}(\phi_2 )(Y_l, Y_r)=\phi_{2_{, y_ly_r}}$, 
\\
{(2.22)}
$\Delta_{2}(\phi_2 )=\sum_{s=1}^{n^2} \epsilon_s \phi_{2_{, y_l y_r}}$ 
\\
{(2.23)}
$|\nabla_{2}(\phi_2 )|^2=\sum_{s=1}^{n_2} \epsilon_s \phi_{2_{, y_s}}^2$.
\\
{(2.24)}
$Hess_{\bar 1}(f_1)(X_i, X_j)=f_{1_{, x_ix_j}}-\sum_k \bar{\Gamma}^k_{ij}f_{1_{, x_k}}$, 
\\
where $\bar{\Gamma}^k_{ij}=0$, $\bar{\Gamma}^i_{ij}=-\frac{\phi_{1_{, x_j}}}{\phi_1}$, $\bar{\Gamma}^k_{ii}=\epsilon_i \epsilon_k\frac{\phi_{1_{, x_k}}}{\phi_1}$ and $\bar{\Gamma}^i_{ii}=-\frac{\phi_{1_{, x_j}}}{\phi_1}$, so (2.24) becomes:
\\
{(2.25)}
$Hess_{\bar 1}(f_1)(X_i, X_j)=f_{1_{, x_ix_j}}+\frac{\phi_{1_{, x_j}}}{\phi_1}f_{1_{, x_i}}+\frac{\phi_{1_{, x_i}}}{\phi_1}f_{1_{, x_j}}$,  for $i \neq j$, and
\\
{(2.26)}
$Hess_{\bar1}(f_1)(X_i, X_i)=f_{1_{, x_ix_i}}+2\frac{\phi_{1_{, x_i}}}{\phi_1}f_{1_{, x_i}}-\epsilon_i \sum_{k=1}^{n_1}\epsilon_k \frac{\phi_{1_{, x_k}}}{\phi_1}f_{1_{, x_k}}$.
\\
\\
Since $Hess_{\bar 2}(f_2)(Y_l, Y_r)=0$, we get:
\\
{(2.27)}
$Hess_{\bar 2}(f_2)(Y_l, Y_r)=f_{2_{, y_ly_r}}+\frac{\phi_{2_{, y_r}}}{\phi_2 }f_{2_{, y_l}}+\frac{\phi_{2_{, y_l}}}{\phi_2 } f_{2_{, y_r}}=0$,  for $l \neq r$, and
\\
{(2.28)}
$Hess_{\bar 2}(f_2)(Y_l, Y_l)=f_{2_{, y_ly_l}}+2\frac{\phi_{2_{, y_l}}}{\phi_2}f_{2_{, y_l}}-\epsilon_l \sum_{s=1}^{n_2}\epsilon_s \frac{\phi_{2_{, y_s}}}{\phi_2}f_{2_{, y_s}}=0$.
\\
\\
Then the Ricci tensors are: 
\\
{(2.29)}
$Ric_{\bar B_1}(X_i, X_j)=\frac{(n_1-2)\phi_{1_{, x_ix_j }}}{\phi_1}$, for $i \neq j$,
\\
{(2.30)}
$Ric_{\bar B_1}(X_i, X_i)=\frac{(n_1-2)\phi_{1_{, x_ix_i }}+ \epsilon_i \sum_{k=1}^{n_1}\epsilon_k\phi_{1_{, x_kx_k}}}{\phi_1}-(n_1-1)\epsilon_i \sum_{k=1}^{n_1}\frac{\epsilon_k \phi_{1_{, x_k}}^2}{\phi_1^2}$,
\\
{(2.31)}
$Ric_{\bar B_2}(Y_l, Y_r)=\frac{(n_2-2)\phi_{2_{, y_ly_r }}}{\phi_2}$, for $l \neq r$,
\\
{(2.32)}
$Ric_{\bar B_2}(Y_l, Y_l)=\frac{(n_2-2)\phi_{2_{, y_ly_l}}+ \epsilon_l \sum_{s=1}^{n_2}\epsilon_s \phi_{2_{, y_sy_s}}}{\phi_2 }-(n_2-1)\epsilon_l \sum_{s=1}^{n_2}\frac{\epsilon_s \phi_{2_{, y_s}}^2}{\phi_2^2}$.
\\
Using (2.29) and (2.25) in the (2.12) and then using (2.30) and (2.26) in the (2.12) we obtain respectively:
\\
\numberwithin{equation}{section}
{(2.33)}
$Ric_{M}(X_i, X_j)=\frac{(n_1-2)\phi_{1_{, x_ix_j }}}{\phi_1} - \frac{d}{f}[f_{1_{, x_ix_j}}+\frac{\phi_{1_{, x_j}}}{\phi_1}f_{1_{, x_i}}+\frac{\phi_{1_{, x_i}}}{\phi_1}f_{1_{, x_j}}]$,  for $i \neq j$,
\\
\numberwithin{equation}{section}
{(2.34)}
$Ric_{M}(X_i, X_i)=\frac{(n_1-2)\phi_{1_{, x_ix_i }}+ \epsilon_i \sum_{k=1}^{n_1}\epsilon_k\phi_{1_{, x_kx_k}}}{\phi_1}-(n_1-1)\epsilon_i \sum_{k=1}^{n_1}\frac{\epsilon_k \phi_{1_{, x_k}}^2}{\phi_1^2} +$
\\
{\centerline{$- \frac{d}{f}[f_{1_{, x_ix_i}}+2\frac{\phi_{1_{, x_i}}}{\phi_1}f_{1_{, x_i}}-\epsilon_i \sum_{k=1}^{n_1}\epsilon_k \frac{\phi_{1_{, x_k}}}{\phi_1}f_{1_{, x_k}}],$}}
\\
while, using (2.31) and (2.27) in the (2.13) and then using (2.32) and (2.28) in the (2.13) we obtain respectively:
\\
\numberwithin{equation}{section}
{(2.35)}
$Ric_{M}(Y_l, Y_r)=\frac{(n_2-2)\phi_{2_{, y_ly_r }}}{\phi_2}$,  for $l \neq r$,
\\
\numberwithin{equation}{section}
{(2.36)}
$Ric_{M}(Y_l, Y_l)=\frac{(n_2-2)\phi_{2_{, y_ly_l}}+ \epsilon_l \sum_{s=1}^{n_2}\epsilon_s \phi_{2_{, y_sy_s}}}{\phi_2 }-(n_2-1)\epsilon_l \sum_{s=1}^{n_2}\frac{\epsilon_s \phi_{2_{, y_s}}}{\phi_2^2}$.
\\
\\
Now considering:
\\
\numberwithin{equation}{section}
{(2.37)}
$Ric_{F}=0$,
\\
\numberwithin{equation}{section}
{(2.38)}
$g_M(U_i, U_j)=f^2g_F(U_i, U_j)$, with $f=f_1+ f_2$,
\\
\numberwithin{equation}{section}
{(2.39)}
$\Delta_{\bar 2}(f_2 )=0$
\\
\numberwithin{equation}{section}
{(2.40)}
$\Delta_{\bar 1}(f_1)=\phi_1^2\sum_{k=1}^{n_1} \epsilon_k f_{1_{, x_k x_k}}-(n_1-2)\phi_1\sum_{k=1}^{n_1} \epsilon_k \phi_{1_{, x_k}} f_{1_{, x_k}}$, 
\\
\numberwithin{equation}{section}
{(2.41)}
 $g_M(\nabla f, \nabla f)=\phi_1^2 \sum_{k=1}^{n_1}\epsilon_k f_{1_{ ,x_k}}^2 + \phi_2^2 \sum_{s=1}^{n_2}\epsilon_s f_{2_{ ,y_s}}^2$,
\\
and by replacing them in (2.17): 
\\
\numberwithin{equation}{section}
{(2.42)}
$Ric_{M}(U_i, U_j)= \{ -f \phi_1^2\sum_{k=1}^{n_1} \epsilon_k f_{1_{, x_k x_k}}+(n_1-2) f \phi_1\sum_{k=1}^{n_1} \epsilon_k \phi_{1_{, x_k}} f_{1_{, x_k}} +$ 
\\
{\centerline {$-(d-1)(\phi_1^2 \sum_{k=1}^{n_1}\epsilon_k f_{1_{ ,x_k}}^2 + \phi_2^2 \sum_{s=1}^{n_2}\epsilon_s f_{2_{ ,y_s}}^2)\}g_F(U_i, U_j)$.}}
\\
\\
Using the equations (2.33), (2.34), (2.35), (2.36) and (2.42), it follows that $(M, g_M)$ is an Einstein manifold if and only if, the equations \textit{(I), (II), (III), (IV), (V)} are satisfied.$\qed$
\\
\\
{\bfseries \centerline {3. The positive constant Ricci curvature case ($\lambda>0$) }}
\\
\\
In this section we look for the existence of a solution to the positive constant Ricci curvature case ($\lambda>0$) when the base-manifold is the product of two conformal manifolds to a $n_1$-dimensional and $n_2$-dimensional pseudo-Euclidean space, respectively, invariant under the action of a $(n_1-1)$-dimensional group of transformations and that the fiber $F$ is flat. 
\\
\\
{\bfseries Theorem 3.1:} \textit{Let $(B, g_B)$, be the base-manifold $B=( B_1 \times B_2)$, $B_1=\mathbb{R}^{n_1}$, with coordinates $(x_1, x_2, .. x_{n_1})$, $B_2=\mathbb{R}^{n_2}$, with coordinates $(y_1, y_2, .. y_{n_2})$, where  $n_1, n_2 \geq 3$, and let  $g_B= g_{B_1} + g_{B_2}$ be the metrics on $B$, where $g_{B_1}=\epsilon_i \delta_{ij}$ and $g_{B_2}=\epsilon_l \delta_{lr}$. 
\\
Let $f_1:\mathbb{R}^{n_1} \rightarrow \mathbb{R}$, $f_2: \mathbb{R}^{n_2} \rightarrow \mathbb{R}$, $\phi_1:\mathbb{R}^{n_1} \rightarrow \mathbb{R}$ and ${\phi_2}:\mathbb{R}^{n_2} \rightarrow \mathbb{R}$, be smooth functions $f_1(\xi_1)$, $f_2(\xi_2)$, $\phi_1(\xi_2)$ and $\phi_2(\xi_2)$,  such that $f(\xi_1, \xi_2)=f_1(\xi_1) + f_2(\xi_2)$ be as in \textit{Definition 1.1}, where $\xi_1=\sum_{i=1}^{n_1} \alpha_i x_i$, $\alpha_i \in \mathbb{R}$, and $\sum_i \epsilon_i \alpha^2_i=\epsilon_{i0}$ or $\sum_i \epsilon_i \alpha^2_i=0$, and by the same token $\xi_2=\sum_{l=1}^{n_2} \alpha_l y_l$, $\alpha_l \in \mathbb{R}$, and $\sum_l \epsilon_l \alpha^2_l=\epsilon_{l0}$ or $\sum_l \epsilon_l \alpha^2_l=0$. 
\\
Finally, let $(M, g_M)$ be $((B_1 \times B_2)\times_{f=f_1+ f_2}F, g_M)$, with $g_M=\bar g_B+(f_1+ f_2)^2 g_F$, with conformal metric $\bar g_B=\bar g_{B_1} + \bar g_{B_2}$, where $\bar g_{B_1}=\frac{1}{\phi_1^2}g_{B_1}$, $\bar g_{B_2}=\frac{1}{{\phi_2}^2}g_{B_2}$, and $F=\mathbb{R}^{d}$ with $g_F=-\delta_{ab}$.
\\
Then, whenever $\sum_i \epsilon_i \alpha^2_i=\epsilon_{i0}$ (and $\sum_l \epsilon_l \alpha^2_l=\epsilon_{l0}$), the warped-product metric \\ $g_M=\bar g_B+(f_1+f_2)^2 g_F$ is Einstein with constant Ricci curvature $\lambda$ if and only if the functions $f_1$, $f_2$, $\phi_1$ and $\phi_2$ satisfy the following conditions:
\\
\\
\textit{(Ia)} $(n_1-2)f\phi_1''-\phi_1 f_1''d - 2\phi_1'f_1'd=0$, for $i \neq j$,
\\
\\
\textit{(IIa)} $\phi_2''=0$, for $l \neq r$,
\\
\textit{(IIIa)} $\sum_k \epsilon_k \alpha^2_k [f \phi_1\phi_1'' - (n_1-1)f \phi_1'^2 + \phi_1\phi_1'f_1'd]=\lambda f$,
\\
\\
\textit{(IVa)} $\sum_s \epsilon_s \alpha^2_s [- (n_2-1) \phi_2'^2]=\lambda$
\\
\\
\textit{(Va)} $\sum_k \epsilon_k \alpha^2_k [-f \phi_1^2f_1''+(n_1-2)f \phi_1\phi_1'f_1'-(d-1)\phi_1^2f_1'^2]+$
\\
{\centerline{$ -\sum_s \epsilon_s \alpha^2_s[(d-1)\phi_2^2 f_2'^2]=\lambda f^2$.}}}
\\
\\
\textit{Proof.} We have:
\\
$\phi_{1_{, x_i x_j}}=\phi_1''\alpha_i \alpha_j$, \quad
$\phi_{1_{, x_i}}= \phi_1' \alpha_i$,  \quad
$ f_{1_{, x_i x_j}}= f_1'' \alpha_i \alpha_j$,  \quad
$f_{1_{, x_i}}=f_1'\alpha_i$, 
\\[0.2cm]
and
\\[0.2cm]
$\phi_{2_{, y_ly_r}}=\phi_2''\alpha_l \alpha_r$,  \quad
$\phi_{2_{, y_l}}= \phi_2' \alpha_l$, \quad
$f_{2_{, y_l y_r}}=f_2'' \alpha_l \alpha_r$, \quad
$f_{2_{, y_l}}=f_2'\alpha_l$.
\\
\\
Substituting these in \textit{(I)} and \textit{(II)} and if $i \neq j$ and $l \neq r$ such that $\alpha_i \alpha_j \neq 0$ and $\alpha_l \alpha_r \neq 0$, we obtain  \textit{(Ia)} and \textit{(IIa)}.
\\
In the same manner for \textit{(III)} and \textit{(IV)}, by considering the relation between $\phi_1''$ and $f_1''$ from \textit{(Ia)} and $\phi_2 ''=0$ from \textit{(IIa)}, we get \textit{(IIIa)} and \textit{(IVa)} respectively. Analogously, the equation \textit{(V)} reduces to \textit{(Va)}.  \qed
\\
\\
Now we are going to look for the existence of a solution to the positive constant Ricci curvature case ($\lambda>0$), considering $f_2(\xi_2)= 1$, and $dim(B_1)=dim(F)$, i.e., $n_1=d$ . So, whenever $\sum_{i=1}^{n_1}\alpha_i^2\epsilon_i \neq 0$, without loss of generality, we may consider $\sum_{i=1}^{n_1}\alpha_i^2\epsilon_i=-1$ (the same for $\sum_{l=1}^{n_2}\alpha_l^2\epsilon_l \neq 0$, in which we consider $\sum_{l=1}^{n_2}\alpha_l^2\epsilon_l=-1$).
\\
In this way the equations \textit{(Ia)}, \textit{(IIa)}, \textit{(IIIa)}, \textit{(IVa)} \textit{(Va)} become:
\\
\\
\textit{(Ib)} $(n_1-2)(f_1+1)\phi_1''-n_1\phi_1f_1''-2n_1\phi_1'f_1'=0$, for $i \neq j$,
\\
\\
\textit{(IIb)}  $\phi_2''=0$, for $l \neq r$,
\\
\\
\textit{(IIIb)} $-(f_1+1) \phi_1\phi_1'' + (n_1-1)(f_1+1) \phi_1'^2 - n_1\phi_1\phi_1' f_1'=\lambda (f_1+1)$,
\\
\\
\textit{(IVb)} $(n_2 -1)\phi_2'^2= \lambda$,
\\
\\
\textit{(Vb)} $(f_1+1)\phi_1^2 f_1''-(n_1-2)(f_1+1)\phi_1 \phi_1'f_1' + (n_1-1)\phi_1^2 f_1'^2=\lambda (f_1+1)^2$.
\\
\\
Note that since $f_2(\xi_2)= constant$, then the equations (2.27) and (2.28), concerning the condition $Hess_{\bar 2}(f_2)=0$, are obviously satisfied. 
\\
It is worth noticing that there is no reason to believe that any nontrivial solutions exist, since the system is overdetermined. One must first check out the compatibility conditions and fortunately this is easy to figure out. Changing the notation:  from $\bigl( \xi_1,\phi_1(\xi_1),f_1(\xi_1)\bigr)$, to $\bigl(t, \beta(t), \gamma(t){-}1\bigr)$ (in order to simplify the writing and avoid confusion with the indexes), and also writing $\lambda = qm^2/2>0$, where $q=n_1$, i.e. $dim(B_1)$, our system of equations then becomes:
\\
\\
\numberwithin{equation}{section}
{(3.1)}$
\begin{cases}
(q-2)\gamma \beta''-q \beta \gamma'' - 2q \beta' \gamma'=0
\\
-\beta \gamma \beta'' - (q- 1) \gamma {\beta'}^2-q \beta'  \gamma' -\tfrac12 qm^2 \gamma=0
\\
\gamma \beta^2  \gamma'' - (q- 2) \beta \gamma \beta'  \gamma' + (q- 1) \beta^2{\gamma'}^2-\tfrac12 qm^2 \gamma^2=0
\end{cases}$
\\
\\
So, if we solve the second and third equations for $\beta''$ and $\gamma''$ and substituting them into the first equation, we note that the first equation can be replaced by a first order equation, that is:
\\
\\
\numberwithin{equation}{section}
{(3.2)}
$(q{-}2) \gamma^2 {\beta'}^2 - 2q \beta \gamma \beta' \gamma' + q \beta^2 {\gamma'}^2 - qm^2 \gamma^2 =: Z(\beta,\gamma,\beta',\gamma')=0$.
\\
\\
Now, differentiating $Z$ with respect to $t$ and then eliminating $\beta''$ and $\gamma''$ using the second and third equations of (3.1), the resulting expression in $(\beta,\gamma, \beta', \gamma')$ is a multiple of $Z(\beta, \gamma, \beta',\gamma')$. This shows us that the combined system of equations (3.1) and (3.2) satisfies the compatibility conditions, so that the system has solutions, specifically, a $3$-parameter family of them.
\\
If we want to describe these solutions more explicitly, we must note that the equations are $t$-autonomous and have a $2$-parameter family of scaling symmetries. In particular, the equations are invariant under the $3$-parameter group of transformations of the form:
\\
\\
\numberwithin{equation}{section}
{(3.3)}
$\Phi_{a,b,c}(t,\beta,\gamma) = (at{+}c, a\beta, b\gamma)$
\\
\\
where $a$ and $b$ are nonzero constants and $c$ is any constant. In fact, the equation (3.2) implies that there is a function $\omega(t)$ such that 
\\
\\
\numberwithin{equation}{section}
{(3.4)}
$\begin{cases}
 \beta' = \frac{2mq \omega (\omega-1)}{\bigl((q-2)\omega^2-2q \omega + q\bigr)}
\\
\gamma' = \frac{m \gamma \bigl((q-2)\omega^2-q\bigr)}{\beta \bigl((q-2)\omega^2-2q \omega + q\bigr)}
\end{cases}$
\\
\\
and then the second and third equations of (3.1) imply that $\omega$ must satisfy
\\
\numberwithin{equation}{section}
{(3.5)}
$\omega' = \frac{m\bigl(q+2q \omega - (3q-2)\omega^2\bigr)}{\beta}$. 
\\
\\
Conversely, the combined system of (3.4) and (3.5) gives the general solution of the original system. This latter system is easily integrated by the usual separation of variables method, i.e., by eliminating $t$ yields a system of the form: 
\\
\\
\numberwithin{equation}{section}
{(3.6)}
$\frac{d\beta}{\beta} = R(\omega)d\omega$
\\
and
\\
\numberwithin{equation}{section}
{(3.7)}
$\frac{d\gamma}{\gamma} = S(\omega) d\omega$
\\
\\
where $R(\omega)$ and $S(\omega)$ are rational functions of $\omega$.  Writing $\beta$ and $\gamma$ as elementary functions of $\omega$, then we can also write:
\\
\\
\numberwithin{equation}{section}
{(3.8)}
$dt = \beta T(\omega) d\omega$,
\\
\\
where $T$ is a rational function of $\omega$, so that $t$ can be written as a function of $\omega$ by quadrature.  Thus, we have the integral curves in $(t,\beta,\gamma,\omega)$-space in terms of explicit functions.
\\
\\
In conclusion (because of the $3$-parameter family of equivalences of solutions), we can say that in certain sense, these solutions are all equivalent to a finite number of possibilities.
\\
\\
{\bfseries{Remarks:}} As is well known, an Einstein warped product manifold with Riemannian-metric and Ricci-flat fiber-manifold can only admit zero or negative Ricci tensor, $Ric \le 0$.
Here we have shown, that a simple pseudo-Riemannian metric construction allows, an Einstein warped product manifold with Ricci-flat fiber-manifold, to obtain $Ric>0$, and this may find interest, for example, in how to build warped-product spacetime models, with positive curvature, whose fiber is Ricci-flat.
\\
\\
\\
{\bfseries{Data availability Statement:}} Not applicable.

\par \bigskip

{\setstretch{1.15}
\begin{itemize}
\setlength\itemsep{0.65em}
\item $[1]$: pigazzini@topositus.com .
\item $[2]$: cozel@kau.edu.sa .
\item $[3]$: saeidjafari@topositus.com .
\item $[4]$: pincak@saske.sk .
\item $[5]$: adebened@sfu.ca .
\end{itemize}
}

\end{document}